\documentclass[hidelinks, 11pt]{article}
\usepackage{amsmath, amssymb, amsfonts, amsthm, mathrsfs, relsize, mathtools, graphicx, float, hyperref, yhmath}
\usepackage[labelsep=period, font=footnotesize, labelfont=bf]{caption}
\usepackage{titlesec}
\titlelabel{\thetitle.\,}

\theoremstyle{plain}
\newtheorem{theorem}{Theorem}[section]
\newtheorem{lemma}[theorem]{Lemma}
\newtheorem{corollary}[theorem]{Corollary}

\theoremstyle{definition}
\newtheorem{remark}[theorem]{Remark}

\DeclareMathAlphabet{\mathbbmsl}{U}{bbm}{m}{sl}
\newcommand{\bmi}[1]{\mbox{\boldmath $ #1$}}
\DeclareMathAlphabet{\mathsfit}{T1}{\sfdefault}{\mddefault}{\sldefault}\SetMathAlphabet{\mathsfit}{bold}{T1}{\sfdefault}{\bfdefault}{\sldefault}

\makeatletter
\newcommand\appendix@section[1]{%
\refstepcounter{section}%
\orig@section*{Appendix \@Alph\c@section. #1}%
\addcontentsline{toc}{section}{Appendix \@Alph\c@section. #1}%
}
\let\orig@section\section
\g@addto@macro\appendix{\let\section\appendix@section}
\makeatother

\makeatletter
\DeclareRobustCommand\widecheck[1]{{\mathpalette\@widecheck{#1}}}
\def\@widecheck#1#2{%
\setbox\z@\hbox{\m@th$#1#2$}%
\setbox\tw@\hbox{\m@th$#1%
\widehat{%
\vrule\@width\z@\@height\ht\z@
\vrule\@height\z@\@width\wd\z@}$}%
\dp\tw@-\ht\z@
\@tempdima\ht\z@ \advance\@tempdima2\ht\tw@ \divide\@tempdima\thr@@
\setbox\tw@\hbox{%
\raise\@tempdima\hbox{\scalebox{1}[-1]{\lower\@tempdima\box
\tw@}}}%
{\ooalign{\box\tw@ \cr \box\z@}}}
\makeatother

\title{\bf{Threshold for  stability of weak saturation}\\[4mm]}

\author{
{\normalsize Mohammadreza Bidgoli}$^{^1}$  \, \,    {\normalsize Ali Mohammadian}$^{^2}$   \, \,   {\normalsize   Behruz Tayfeh-Rezaie}$^{^1}$  \, \,  {\normalsize  Maksim Zhukovskii}$^{^{3,4,5}}$ \\[6mm]
$^{^1}${\normalsize School of Mathematics,} \\
{\normalsize Institute for Research in Fundamental Sciences\,(IPM),} \\
{\normalsize P.O. Box 19395-5746, Tehran, Iran}\\[2mm]
$^{^2}${\normalsize School of Mathematical Sciences, Anhui University,} \\
{\normalsize Hefei 230601, Anhui,  China} \\[2mm]
$^{^3}${\normalsize Laboratory of Combinatorial and Geometric Structures, MIPT,} \\
{\normalsize Moscow Region 141701, Russian Federation} \\[2mm]
$^{^4}${\normalsize Adyghe State University, Caucasus Mathematical Center,} \\
{\normalsize Republic of Adygea, 385000, Russian Federation} \\[2mm]
$^{^5}${\normalsize RANEPA, Moscow, 119571, Russian Federation} \\[4mm]
{\normalsize \texttt{bd@ipm.ir}}  \, \,  {\normalsize \texttt{ali\_m@ahu.edu.cn}} \, \,   {\normalsize \texttt{tayfeh-r@ipm.ir}}  \, \,    {\normalsize \texttt{zhukmax@gmail.com}}\\[6mm]
}

\date{}

\begin{document}

\maketitle

\begin{abstract}
We study  the weak $K_s$-saturation number of the Erd\H{o}s--R\'{e}nyi random graph $\mathbbmsl{G}(n, p)$, denoted by $\mathrm{wsat}(\mathbbmsl{G}(n, p), K_s)$, where $K_s$ is the complete graph on $s$ vertices. Kor\'{a}ndi and Sudakov  in 2017   proved that the weak $K_s$-saturation number of $K_n$ is stable,  in the sense that    it remains the same after removing edges with constant probability. In this paper, we prove that there exists a threshold  for this stability property and give upper and lower bounds on the threshold.  This generalizes    the result  of Kor\'{a}ndi  and Sudakov.   A general upper bound for  $\mathrm{wsat}(\mathbbmsl{G}(n, p), K_s)$ is also provided. \\[3mm]
\noindent {\bf AMS Mathematics Subject Classification\,(2020):}    05C80, 05C35, 60K35. \\[6mm]
\end{abstract}

\section{Introduction}

Given   a graph $F$, an {\sl $F$-bootstrap percolation process} is a sequence of graphs $H_0\subset H_1\subset\cdots\subset H_m$ such that, for   $i=1,\ldots,m$, $H_i$ is obtained from $H_{i-1}$ by adding an edge that belongs to a copy of $F$ in $H_i$. The $F$-bootstrap percolation process was introduced by Bollob\'{a}s more than 50 years ago \cite{bol} and can be seen as a special case of the `cellular automata' introduced by von Neumann \cite{Neumann} after a suggestion of Ulam \cite{Ulam}. The  $F$-bootstrap percolation is also similar to $r$-neighborhood bootstrap percolation model having applications in physics; see, for example, \cite{Adler}, \cite{Fontes}   and \cite{Morris}.

Given two graphs $G$ and $F$,  a spanning subgraph  $H$  of $G$  is  said to be a  {\sl weakly $F$-saturated subgraph of $G$} if $H$ contains no subgraph  isomorphic to $F$ and there exists an $F$-bootstrap percolation process $H=H_0\subset H_1\subset\cdots\subset H_m=G$. The  minimum number of edges  in  a weakly $F$-saturated  subgraph of $G$ is called the {\sl weak $F$-saturation number} of $G$ and is denoted by $\mathrm{wsat}(G, F)$.

We denote by $\mathbbmsl{G}(n,p)$  the Erd\H{o}s--R\'{e}nyi  random graph on vertex set $[\![n]\!]=\{1,\ldots,n\}$ constructed by adding every edge $e\in\{xy \, | \, x, y\in [\![n]\!] \text{ and }  x\neq y\}$
with probability    $p$  independently of all the others. Kor\'{a}ndi and   Sudakov    \cite{kor}  initiated  the study of  weak saturation numbers of random graphs. They   proved that, for every fixed real number  $p\in(0, 1)$ and integer  $s\geq 3$, $\mathrm{wsat}(\mathbbmsl{G}(n, p), K_s)=\mathrm{wsat}(K_n, K_s)$ with high probability. Recall that the  notion `with high probability', which is written as  `whp' for brevity,   is used whenever   an event  occurs in     $\mathbbmsl{G}(n,p)$  with a  probability approaching $1$ as $n\to\infty$. It was already known that    $\mathrm{wsat}(K_n, K_s)={s-2 \choose 2}+(s-2)(n-s+2)$ by a classic result proved by Lov\'{a}sz \cite{L77}.
Other proofs of this result have been obtained by Frankl \cite{Frankl82}, Kalai \cite{Kalai84, Kalai85}, Alon \cite{Alon85}  and Yu \cite{Yu93}.
Kor\'{a}ndi and Sudakov    \cite{kor} also     noticed that  $\mathrm{wsat}(\mathbbmsl{G}(n, p), K_s)=\mathrm{wsat}(K_n, K_s)$ whp when $p\geq n^{-\varepsilon}$ for small
enough $\varepsilon> 0$ and asked about smaller $p$ and about possible  threshold probability for the
property of having the weak $K_s$-saturation number of $\mathbbmsl{G}(n,p)$  exactly ${s-2 \choose 2}+(s-2)(n-s+2)$.    We denote this property   by $\mathcal{A}_s$.
In this paper, we prove that this threshold exists and present  upper and lower  bounds on that. The formal   definition of  a    threshold function      appears  in Page 18 of    \cite{RG}.

The rest of the paper is organized as follows. In Section \ref{prelim}, we fix some  notation     used in the paper and state  the known results   that we have referred to. In Section \ref{secLowerP},  we prove that there is a threshold probability  for the property   $\mathcal{A}_s$ and present a lower
bound on it. This is done by considering  some  auxiliary events. We establish an upper bound for the aforementioned threshold in Section \ref{secUpperP} by introducing a weakly $K_s$-saturated subgraph of $\mathbbmsl{G}(n,p)$. Finally, in Section \ref{secUpperwsat}, we find a universal upper bound for $\mathrm{wsat}(\mathbbmsl{G}(n, p), K_s)$ to also cover the gap of $p$ between the provided upper and lower bounds on the threshold.

\section{Preliminaries}\label{prelim}

In this section, we introduce notation  and formulate several probabilistic inequalities  that we use in the rest of the  paper.

For a graph $G$, we denote  the vertex set and the edge set of   $G$  by $V(G)$ and $E(G)$, respectively.
The  {\sl size} of $G$ is  defined as   $|E(G)|$ and is denoted by $e(G)$.
For a vertex $v$ of $G$, set $N_G(v)=\{x\in V(G) \, | \, v \text{  is adjacent to } x\}$ to be the neighborhood of   $v$ in $G$. Also, for a subset $U$ of $V(G)$, define the open and closed neighborhood of $U$ in $G$ as $N_G(U)=\cap_{u\in U} N_G(u)$  and $N_G[U]=U \cup N_G(U)$, respectively.
Furthermore, for  a subset $S$ of $V(G)$, we denote the induced subgraph of $G$ on $S$  by $G[S]$.

We also use the standard asymptotic notation in the rest of the paper. For two real-valued functions $f(n)$ and $g(n)$, we write $f(n)=O(g(n))$ if there exists a constant $c>0$   such that $|f(n)|\leq cg(n)$ for every   large enough integer $n$. Also, we use the notation  $f(n)=o(g(n))$ if the same holds for any constant $c>0$.    We sometimes  write    $f(n)\ll g(n)$  and   $g(n)=\omega(f(n))$ instead of   $f(n)=o(g(n))$. Finally, we use the notation   $f(n)=\mathnormal{\Theta}(g(n))$ if both   $f(n)=O(g(n))$ and $g(n)=O(f(n))$  hold.

In what follows, we   recall the  probabilistic inequalities that we make use of all in the  next sections.

\begin{theorem}[Markov's inequality;  Inequality (1.3)  of  \cite{RG}]\label{Markov}
Let   $X$ be a  nonnegative  random variable. Then,  for all  $t>0$,
$$\mathbbmsl{P}[X\geq t]\leq\frac{\mathbbmsl{E}[X]}{t}.$$
\end{theorem}

\begin{corollary}\label{First-Moment}
Let   $X$ be a   nonnegative  integer-valued       random variable. If   $\mathbbmsl{E}[X]=o(1)$, then $X=0$  whp.
\end{corollary}

\begin{theorem}[Chebyshev's inequality;  Inequality (1.2)  of  \cite{RG}]\label{chebyshev}
Let   $X$ be a  random variable with the expected value    $\mathbbmsl{E}[X]$ and   the variance  $\mathbbmsl{Var}[X]$. Then,  for all  $t>0$,
$$\mathbbmsl{P}\Big[\big|X-\mathbbmsl{E}[X]\big|\geq t\Big]\leq\frac{\mathbbmsl{Var}[X]}{t^2}.$$
\end{theorem}

\begin{corollary}\label{Second-Moment}
Let   $X$ be a     random variable   with positive   expected value. If   $\mathbbmsl{Var}[X]=o(\mathbbmsl{E}[X]^2)$, then $X>0$  whp.
\end{corollary}

\begin{theorem}[Chernoff's inequality;  Theorem 2.1 of  \cite{RG}]\label{chernoff}
Let   $X\sim\mathrm{Bin}(n,p)$ be a binomial random variable  with parameters $n$ and $p$. Then, for any  $t\geq 0$,
$$\mathbbmsl{P}[X\leq np-t]\leq\exp \left(-\frac{t^2}{2np}\right).$$
\end{theorem}

The following consequence  of the  Fortuin--Kasteleyn--Ginibre  inequality  \cite[Theorem 2.12]{RG}  appears in  Page 31 of \cite{RG}.

\begin{theorem}\label{FKG}
Let $\mathcal{S}$ be a family of subgraphs of $K_n$ and assume that   the random variable $X$    counts  the number of graphs in  $\mathcal{S}$ that appear in $\mathbbmsl{G}(n,p)$. Then,
$$\mathbbmsl{P}[X=0]\geq\prod_{H\in \mathcal{S}}\Big(1-\mathbbmsl{P}\big[H\subset\mathbbmsl{G}(n,p)\big]\Big).$$
\end{theorem}

\begin{theorem}[Janson's inequality;  Theorem 2.18 of \cite{RG}]\label{Janson}
Let $\mathcal{S}$ be a family of subgraphs of $K_n$. Assume that   the random variable $X$    counts  the number of graphs in  $\mathcal{S}$ that appear in $\mathbbmsl{G}(n,p)$. For every $H_1,H_2\in\mathcal{S}$, let $H_1\sim H_2$ indicate  that $H_1\neq H_2$ and $H_1, H_2$  share at least one edge. Define
$$\mathnormal{\Delta}=\mathop{\sum_{H_1, H_2\in\mathcal{S}}}_{H_1\sim H_2}\mathbbmsl{P}\big[H_1, H_2\subset\mathbbmsl{G}(n,p)\big].$$
Then,
$$\mathbbmsl{P}[X=0]\leq\exp\left(-\mathbbmsl{E}[X]+\frac{\mathnormal{\Delta}}{2}\right).$$
\end{theorem}

\section{The existence of the   threshold}\label{secLowerP}

In this section, we prove the existence of the  threshold probability  for the property   $\mathcal{A}_s$  and present a lower bound on it.

\begin{theorem}\label{lower-ii}
For a  fixed integer  $s\geq 3$,  let
$$c_s=\left(2\left(1-\frac{1}{s-2}\right)(s+1)!\right)^{\frac{2}{(s-2)(s+1)}}$$  and $$q_s(n)=n^{-\frac2{s+1}}(\ln n)^\frac{2}{(s-2)(s+1)}.$$
If $p\leq c_sq_s$, then   the property   $\mathcal{A}_s$ does not hold in $\mathbbmsl{G}(n, p)$ whp.
\end{theorem}

\begin{proof}
Let $X_s$ be the number of $K_s$  in $\mathbbmsl{G}(n,p)$.
If $p\leq\frac{1}{n\ln n}$, then    $e(\mathbbmsl{G}(n, p))<\frac{n}{\sqrt{\ln n}}$  whp  using  Theorem \ref{Markov}.  As  $\mathrm{wsat}(\mathbbmsl{G}(n, p), K_s)\leq e(\mathbbmsl{G}(n, p))$, we get that $\mathrm{wsat}(\mathbbmsl{G}(n, p), K_s)\neq\mathrm{wsat}(K_n, K_s)$ whp.
Now, let $\frac{1}{n\ln n}<p\ll n^{-2/(s+1)}$ and let $\varepsilon\in(0, 1)$ be a small constant. By Theorem \ref{chebyshev}, the number of edges of $\mathbbmsl{G}(n,p)$ belongs to $(\tfrac{n^2p}{2+2\varepsilon}, \tfrac{n^2p}{2-2\varepsilon})$  whp,   since it has binomial distribution with parameters ${n\choose 2}$ and $p$. Consider an arbitrary increasing sequence $w_n=\omega(1)$ such that $(np^{(s+1)/2})^{s-2}w_n=o(1)$. Using  Theorem \ref{Markov},
$$\mathbbmsl{P}\left[X_s>\frac{n^2p}{w_n}\right]\leq\frac{\mathbbmsl{E}[X_s]}{\frac{n^2p}{w_n}}<\frac{n^sp^{s\choose 2}w_n}{n^2p}=\left(np^{\frac{s+1}{2}}\right)^{s-2}w_n\to0,$$
implying that $X_s=o(n^2p)$ whp. Since $e(\mathbbmsl{G}(n, p))-X_s\leq\mathrm{wsat}(\mathbbmsl{G}(n, p), K_s)\leq e(\mathbbmsl{G}(n, p))$, we get   that  $\mathrm{wsat}(\mathbbmsl{G}(n, p), K_s)=e(\mathbbmsl{G}(n, p))(1+o(1))$ whp. From this, we immediately deduce that   $\mathrm{wsat}(\mathbbmsl{G}(n, p), K_s)\neq\mathrm{wsat}(K_n, K_s)$ whp when  $p\ll n^{-2/(s+1)}$ and  $p\notin I$, where $$I=\left(\frac{2(1-\varepsilon)(s-2)}{n}, \frac{2(1+\varepsilon)(s-2)}{n}\right).$$

If  $s\geq4$ and $p\in I$, then  $X_s=0$ whp, since $p<2(1+\varepsilon)(s-2)/n$ and $n^{-2/(s-1)}$ is the threshold for appearance of $K_s$   by  Theorem 3.4 of   \cite{RG}. Hence,    the weak saturation number is exactly $e(\mathbbmsl{G}(n, p))$ that is not concentrated on a single value.

If $s=3$ and $p\in[\tfrac{2-2\varepsilon}{n}, \tfrac{2+2\varepsilon}{n}]$,  then $\mathrm{wsat}(\mathbbmsl{G}(n, p), K_s)$ whp       equals the difference between the number of edges and the number of triangles. Note that whp all the triangles in $\mathbbmsl{G}(n, p)$ are disjoint,  since whp there are no subgraphs with at most $5$ vertices and at least two cycles  using   Theorem 3.4 of   \cite{RG}. The latter random variable is also not concentrated in a unit set. This is because the number of edges has the binomial distribution with parameters ${n\choose 2}$ and $p=\mathnormal{\Theta}(\tfrac{1}{n})$,  so
it is outside any interval of length  $o(\sqrt{n})$ whp,  while the number of triangles is bounded from above by an asymptotically Poisson random variable.  For a constant nonnegative  integer  $L$, the property of having at most $L$ triangles is   decreasing,  implying that $\mathbbmsl{P}[X_3\leq L]$ is minimum when $p=\tfrac{2+2\varepsilon}{n}$. At the same time,  the number of  triangles in $\mathbbmsl{G}(n, \tfrac{2+2\varepsilon}{n})$ converges in distribution to
a Poisson random variable with parameter $\tfrac{4}{3}(1+\varepsilon)^3$
using    Lemma 1.10 and Theorem 3.19 of \cite{RG}. Thus,   for an arbitrary slowly increasing function $z_n=\omega(1)$,   $X_3<z_n$ whp.

Summing up, for $p\ll n^{-2/(s+1)}$,   $\mathbbmsl{G}(n,p)$ does not have $\mathcal{A}_s$ whp.
Now, let $\gamma>0$ be fixed and small enough  and  let $p\geq n^{-2/(s+1)-\gamma}$. Also, denote  by $\mathcal{B}_s$ the property  that every edge belongs to some $K_s$. To continue  the proof, we need the following technical lemma about $\mathcal{B}_s$.

\begin{lemma}\label{edge_in_clique}
Let $w_1, w_2, \ldots$ be a sequence of real numbers  and let
$$p(n)=c_sq_s(n)\left(1+\frac{2\ln\ln n}{s(s-2)^2(s+1)\ln n}+\frac{w_n}{\ln n}\right).$$
Then, the  following hold:
\begin{enumerate}
\item[{\rm (i)}] If $w_n\to\infty$, then   $\mathbbmsl{G}(n,p)$ has $\mathcal{B}_s$ whp.
\item[{\rm (ii)}] If $w_n\to-\infty$, then   $\mathbbmsl{G}(n,p)$ does not have $\mathcal{B}_s$ whp.
\end{enumerate}
\end{lemma}

We include the proof of Lemma \ref{edge_in_clique}  in Appendix   \ref{appen}. To end the proof,  suppose  that  $p\leq c_sq_s$.
We may assume that  $p(n)$ is of the form   given  in Lemma \ref{edge_in_clique} with  $w_n\to-\infty$.  Lemma  \ref{edge_in_clique} implies that  there exists an edge $e$ in $\mathbbmsl{G}(n, p)$ that is not contained in any $K_s$ which   immediately yields  that $e$ should belong to all   weakly $K_s$-saturated subgraphs of  $\mathbbmsl{G}(n, p)$. Let  $F$ be  a weakly $K_s$-saturated subgraph of $\mathbbmsl{G}(n,p)$ of size ${s-2 \choose 2}+(s-2)(n-s+2)$. Then,  $F-e$  is a weakly $K_s$-saturated subgraph of  $\mathbbmsl{G}(n,p)-e$ of size    ${s-2 \choose 2}+(s-2)(n-s+2)-1$. But, this is impossible,  since   $\mathbbmsl{G}(n, p)-e$ is a  weakly $K_s$-saturated subgraph of  $K_n$ whp. To see this,    note that  the threshold probability  for $\mathbbmsl{G}(n, p)$ to be a  weakly $K_s$-saturated subgraph of  $K_n$ is at most $n^{-2/(s+1)-\gamma}$ for small enough $\gamma>0$ by Theorem 1 of  \cite{bal} and the same argument works for $\mathbbmsl{G}(n,p)-e$. Since $p(n)\geq n^{-2/(s+1)-\gamma}$ for small enough $\gamma>0$, the proof of Theorem \ref{lower-ii} is completed.
\end{proof}

\begin{theorem}\label{lower-i}
For any   fixed integer  $s\geq 3$, there exists a threshold probability for the property   $\mathcal{A}_s$.
\end{theorem}

\begin{proof}
Let  $\mathcal{B}^*_s$ be  the   property that every pair of vertices have   $s-2$ common neighbors that induce a clique. This property is   increasing  with  a    sharp threshold probability  $$q^*_s(n)=(2(s-2)!)^{\frac{2}{(s-2)(s+1)}}q_s(n).$$ The above  threshold  function   is  obtained by    the  same way as in the proof of Lemma \ref{edge_in_clique} and  originally comes from   \cite{SP90}. Clearly, $\mathcal{A}_s\cap\mathcal{B}^*_s$ is increasing as well and thus  it has a threshold probability $r_s(n)$ by Theorem 1.24 of  \cite{RG}. Obviously,  $r_s\geq q^*_s$.
In what follows, we prove that   $r_s$ is a threshold for $\mathcal{A}_s$.

If $p\gg r_s$, then  it follows from   Theorem \ref{lower-ii} that   $\mathbbmsl{G}(n,p)$ has the property $\mathcal{A}_s$ whp. It remains to prove the opposite for $p\ll r_s$.
If $p(n)\leq c_sq_s(n)$ for infinitely many $n$, then  it follows from    Theorem \ref{lower-ii}  that  $\mathbbmsl{G}(n, p)$ does not have the property $\mathcal{A}_s$ whp. So, we may assume $p(n)>c_sq_s(n)$ for all large enough integer  $n$.   If $(1+\varepsilon)q^*_s(n)\leq p(n)\ll r_s(n)$  for  some fixed  $\varepsilon>0$ and  infinitely many $n$, then
$$\mathbbmsl{P}\big[\mathbbmsl{G}(n,p)\in\mathcal{A}_s\big]\leq\mathbbmsl{P}\big[\mathbbmsl{G}(n,p)\notin\mathcal{B}^*_s\big]+\mathbbmsl{P}\big[\mathbbmsl{G}(n,p)\in\mathcal{A}_s\cap\mathcal{B}^*_s\big]\to0$$ as $n\to\infty$.
Finally, suppose  by contradiction that there exists a fixed    $\delta>0$ such that
$c_sq_s(n)<p(n)\leq(1+\varepsilon)q^*_s(n)$ and $\mathbbmsl{P}[\mathbbmsl{G}(n, p)\in\mathcal{A}_s]\geq\delta$
for  some fixed  $\varepsilon>0$ and for any $n$ belong to  an infinite  set $\mathnormal{\Lambda}$ of positive integers.  In the rest of the proof,    we always assume that $n\in\mathnormal{\Lambda}$.
Consider an auxiliary function  $$\widecheck{p}(n)=c_sq_s(n)\left(1+\frac{1}{\sqrt{\ln n}\ln\ln n}\right)$$ and let $\widehat{p}(n)=(1+\varepsilon)q^*_s(n)$ for the sake of simplicity.  If $p<\widecheck{p}$, then $\mathbbmsl{G}(n,\widecheck{p})$ can be obtained from $\mathbbmsl{G}(n,p)$ by drawing every missing edge with probability $(\widecheck{p}-p)/(1-p)$. By the same arguments as given in Appendix \ref{appen},
the probability that there exists an edge $xy$ in $\mathbbmsl{G}(n,\widecheck{p})\setminus\mathbbmsl{G}(n,p)$   such that   $N_{\mathbbmsl{G}(n,p)}(\{x, y\})$ does not contain  a clique of size   $s-2$   is at most
$${n\choose 2}\frac{\widecheck{p}-p}{1-p}\exp\left(-{n-2\choose s-2}\big(c_sq_s(n)\big)^{\frac{(s-2)(s+1)}{2}}+o(1)\right)=O\left(\frac{1}{\ln\ln n}\right).$$
Therefore,
$$\mathbbmsl{P}\big[\mathbbmsl{G}(n,\widecheck{p})\in\mathcal{A}_s\big]\geq\mathbbmsl{P}\big[\mathbbmsl{G}(n,p)\in\mathcal{A}_s\big]-O\left(\frac{1}{\ln\ln n}\right)\geq\delta-o(1).$$
So,  we may assume that $p\geq\widecheck{p}$. In the same manner, $\mathbbmsl{G}(n,\widehat{p})$ can be obtained from $\mathbbmsl{G}(n,p)$  by drawing every missing edge with probability $(\widehat{p}-p)/(1-p)$.
Then, the probability that there exists an edge $xy$ in $\mathbbmsl{G}(n,\widehat{p})\setminus\mathbbmsl{G}(n,p)$
such that   $N_{\mathbbmsl{G}(n,p)}(\{x, y\})$ does not contain  a clique of size   $s-2$     is at most
$${n\choose 2}\frac{\widehat{p}-p}{1-p}\exp\left(-{n-2\choose s-2}\widecheck{p}^{\frac{(s-2)(s+1)}{2}}+o(1)\right)=O\left((\ln n)^{\frac{2}{(s-2)(s+1)}}\exp\left(-\frac{\sqrt{\ln n}}{\ln\ln n}\right)\right).$$
Hence,
$$\mathbbmsl{P}\big[\mathbbmsl{G}(n,\widehat{p})\in\mathcal{A}_s\big]\geq\mathbbmsl{P}\big[\mathbbmsl{G}(n,p)\in\mathcal{A}_s\big]-o(1)\geq\delta-o(1),$$
a contradiction.
\end{proof}

\section{An  upper bound on the  threshold}\label{secUpperP}

In this section, we   present an  upper    bound on   the  threshold probability  for the property   $\mathcal{A}_s$. Let us first  recall the following definition.
The {\sl $k$-th power}  of a graph  $\mathnormal{\Gamma}$, denoted by $\mathnormal{\Gamma}^k$,  is the graph  with    vertex set $V(\mathnormal{\Gamma})$  such that two vertices $x, y$ are adjacent in $\mathnormal{\Gamma}^k$ if and only if the distance between  $x, y$  in  $\mathnormal{\Gamma}$  is at most $k$. We need the following result on the threshold probability of the appearance of the $k$-th power of a Hamilton cycle.

\begin{theorem}[\cite{posa, kahn, Riordan}]
If   $p>c\tfrac{\ln n}{n}$ for a sufficiently large   constant  $c$,  then  $\mathbbmsl{G}(n,p)$ contains a Hamilton cycle whp. For every integer $k\geq 2$ and  $p\gg n^{-1/k}$, then $\mathbbmsl{G}(n,p)$ contains the $k$-th power of a Hamilton cycle whp.
\label{th:Hamilton_powers}
\end{theorem}

We  also  define  two graph properties to use later.
We say     $\mathnormal{\Gamma}$ has the property {\sf EXT} if for every subset $S\subseteq V(\mathnormal{\Gamma})$  of size $s$, $N_\mathnormal{\Gamma}(S)$ contains a clique of size  $s-2$.
Moreover, we say   $\mathnormal{\Gamma}$ has the property {\sf HAM} if for every subset $S\subseteq V(\mathnormal{\Gamma})$  of size $s-1$,  $\mathnormal{\Gamma}[N_\mathnormal{\Gamma}(S)]$ contains the $(s-2)$-th power of a Hamilton path.

\begin{lemma}\label{11}
Let   $s\geq 3$   and     $n\geq s-2$. Assume that  both properties  {\sf EXT} and {\sf HAM} hold for a graph $G$ on $n$ vertices. Then,  $\mathrm{wsat}(G,K_s)\leq{s-2 \choose 2}+(s-2)(n-s+2)$.
\end{lemma}

\begin{proof}
If $n\in\{s-2, s-1\}$, then the  result is clearly valid.    Let $n\geq s$ and let $\mathnormal{\Omega}$ be a clique of size $s-2$  in $G$. We define a  spanning   subgraph $H$ of $G$ as follows.  The  graph  $H$ contains all edges of $G$ with endpoints in $\mathnormal{\Omega}$  and  also   all edges of $G$ with endpoints in both  $\mathnormal{\Omega}$ and $N_G(\mathnormal{\Omega})$. We still have to add to $H$  some other edges going outside $N_G[\mathnormal{\Omega}]$. For every $v\in V(G)\setminus N_G[\mathnormal{\Omega}]$, we add    $s-2$ edges of $G$ adjacent to $v$ described below.  Since $G$   satisfies    {\sf HAM},  the graph $H_v=G[N_G(\{v\}\cup\mathnormal{\Omega})]$ contains  the  $(s-2)$-th power of a Hamilton path. Starting from  a beginning vertex, denote the vertices of $H_v$ going in the natural order induced by the Hamilton path  by $x^v_1, \ldots, x^v_{h_v}$, where $h_v=|V(H_v)|$. We add the  edges $vx^v_1,  \ldots,  vx^v_{s-2}$ to $H$ for any $v\in V(G)\setminus N_G[\mathnormal{\Omega}]$. It is easy to see that $H$ is of size ${s-2 \choose 2}+(s-2)(n-s+2)$, so it suffices to prove that $H$ is a weakly $K_s$-saturated subgraph of  $G$.

First, all edges with endpoints in  $N_G(\mathnormal{\Omega})$ can be activated,  since they  belong to a  $K_s$ containing  $\mathnormal{\Omega}$. Next, for each $v \in V(G)\setminus N_G[\mathnormal{\Omega}]$, we may active  the edges $vx^v_{s-1},  \ldots, vx^v_{h_v}$ one by one,  since every such edge belongs to a $K_s$  containing the previous $s-2$ vertices of the $(s-2)$-th power of the Hamilton path. Finally, each edge $xy$ with endpoints in $V(G)\setminus N_G[\mathnormal{\Omega}]$ can be activated. To see this, note that    $N_G(\{x, y\}\cup\mathnormal{\Omega})$ contains a clique of size $s-2$,  say $\mathnormal{\Omega}_{xy}$, since  $G$   satisfies  {\sf EXT}.  It follows from   $\mathnormal{\Omega}_{xy}\subseteq N_G(\mathnormal{\Omega})$ that the edges  with endpoints in  $\mathnormal{\Omega}_{xy}$ are already activated   and  so $xy$ is the last edge of the  $\{x,y\}\cup\mathnormal{\Omega}_{xy}$ of size $s$ and   can be activated  as well.
\end{proof}

\begin{lemma}\label{12}
Let  $s\geq 3$ be  a fixed  integer and $p\geq n^{-1/(2s-3)}\ln n$. Then,  $$\mathbbmsl{P}\big[\text{$\mathbbmsl{G}(n, p)$  has both properties $\mathsf{EXT}$ and $\mathsf{HAM}$}\big]\to1$$ as $n\to\infty$.
\end{lemma}

\begin{proof}
By Theorem 2 of \cite{SP90},  $\mathbbmsl{G}(n,p)$ has the property {\sf EXT} whp provided that    $$p\gg n^{-\frac{2}{3(s-1)}}(\ln n)^{\frac{2}{3(s-2)(s-1)}}.$$
Therefore, if  $p\geq n^{-1/(2s-3)}\ln n$, then $\mathbbmsl{P}[\text{$\mathbbmsl{G}(n, p)$  has the  property $\mathsf{EXT}$}]\to1$ as $n\to\infty$.
To explore the  behavior of  $\mathbbmsl{P}[\text{$\mathbbmsl{G}(n, p)$  has the  property $\mathsf{HAM}$}]$, we use the    fact that, for every   positive integers $k$ and $r$, if    $m$ is sufficiently large and  $q(m)\gg m^{-1/k}(\ln m)^2$,  then $\mathbbmsl{G}(m, q)$ contains the $k$-th power of a Hamilton path    with probability at least $1-m^{-r}$.
To see  this fact, we apply Theorem \ref{th:Hamilton_powers} and  conclude   that
there exists a constant $C>0$ such that if $\widecheck{q}(m)>Cm^{-1/k}\ln m$, then $\mathbbmsl{G}(m, \widecheck{q})$ contains the $k$-th power of a Hamilton path   with probability at least $1-\tfrac{1}{e}$.
To boost this probability to $1-m^{-r}$,  it suffices to apply    the coupling technique  by  taking  the union of $\lceil r\ln m\rceil$ independent copies of $\mathbbmsl{G}(m, \widecheck{q})$ on vertex set $[\![m]\!]$ and letting  $p={\lceil r\ln m\rceil}\widecheck{q}$.

Fix $W\subseteq[\![n]\!]$ of size $s-1$ and set $m=np^{s-1}/2$. Using  Theorem \ref{chernoff},
the probability that   $N_{\mathbbmsl{G}(n,p)}(W)$ does not contain the $(s-2)$-th power of a Hamilton path is at most
$$\mathbbmsl{P}\Big[\big|N_{\mathbbmsl{G}(n,p)}(W)\big|<m\Big]+m^{-3s}\mathbbmsl{P}\Big[\big|N_{\mathbbmsl{G}(n,p)}(W)\big|\geq m\Big]\leq\exp\left(-\frac{m}{4}\right)+n^{-s}\to 0$$ as $n\to 0$.
Now, the result follows from the union bound theorem.
\end{proof}

In Theorem 1.3 of \cite{kor},  Kor\'{a}ndi  and Sudakov proved, if  $p\in(0, 1)$  is a   constant probability  and    $s\geq3$ is   a fixed  integer, then  $\mathrm{wsat}(\mathbbmsl{G}(n,p), K_s)={s-2\choose 2}+(s-2)(n-s+2)$ whp.  The following theorem  generalizes    their  result.

\begin{theorem}\label{th_upper}
Let  $s\geq 3$ be  a fixed  integer and $p\geq n^{-1/(2s-3)}\ln n$. Then,  $\mathbbmsl{G}(n,p)$ has the property $\mathcal{A}_s$ whp.
\end{theorem}

\begin{proof}
Using    Theorem 1 of  \cite{bal},    $\mathbbmsl{G}(n,p)$ is a weakly $K_s$-saturated subgraph of  $K_n$ whp and thus    we get  that $\mathrm{wsat}(\mathbbmsl{G}(n, p), K_s)\geq\mathrm{wsat}(K_n, K_s)={s-2 \choose 2}+(s-2)(n-s+2)$  whp.
It remains to prove that  there exists a weakly $K_s$-saturated subgraph of $\mathbbmsl{G}(n,p)$ of size   ${s-2\choose 2}+(s-2)(n-s+2)$    whp. So,
the result immediately  follows  from Lemmas \ref{11} and \ref{12}.
\end{proof}

\begin{remark}\label{rmrk}
The power of  the logarithm  factor in Theorem \ref{th_upper} is not the best possible.  By Theorem \ref{th:Hamilton_powers}, we may replace the factor $\ln n$ in Theorem \ref{th_upper} with $(\ln n)^{\xi}$, where
$$\xi=\left\{\begin{array}{ll}
\frac{2}{3}  &   \mbox{ if } s=3\mbox{,}\\\vspace{-3mm}\\
\frac{s-2}{2s-3} &  \mbox{ if } s\geq 4\mbox{.}
\end{array}\right.$$
\end{remark}

\section{An  upper bound on $\bmi{\mathrm{wsat}(\mathbbmsl{G}(n, p), K_s)}$}\label{secUpperwsat}

From the previous arguments, we know that  $\mathrm{wsat}(\mathbbmsl{G}(n, p), K_s)=e(\mathbbmsl{G}(n,p))(1+o(1))$ whp when $p(n)\ll n^{-2/(s+1)}$ and $\mathrm{wsat}(\mathbbmsl{G}(n, p), K_s)={s-2\choose 2}+(s-2)(n-s+2)$  whp when $p(n)\geq n^{-1/(2s-3)}(\ln n)^{\xi}$ for    $\xi$   presented  in Remark \ref{rmrk}. In this section, we give  an upper bound on  $\mathrm{wsat}(\mathbbmsl{G}(n, p), K_s)$ for all the remaining values of $p$. It worth to mention  that $\mathrm{wsat}(\mathbbmsl{G}(n, p), K_s)\geq{s-2\choose 2}+(s-2)(n-s+2)$ whp,  since $\mathbbmsl{G}(n,p)$ is a  weakly $K_s$-saturated subgraph of  $K_n$  whp  for those values of $p$ by Theorem 1 of \cite{bal}.

\begin{theorem}\label{upperbd}
Let  $s\geq 3$ be  a fixed  integer and let
$$\gamma=\left\{\begin{array}{ll}
1  &   \mbox{ if } s=3\mbox{,}\\\vspace{-3mm}\\
0 &  \mbox{ if } s\geq 4\mbox{.}
\end{array}\right.$$
Also, let  $w_1, w_2, \ldots$ be a sequence of real numbers such that   $w_n\to\infty$ as $n\to\infty$. Then,      $$\mathrm{wsat}\big(\mathbbmsl{G}(n, p), K_s\big)\leq (s-2)n+\frac{w_n(\ln n)^{2(\gamma+s-2)}}{p^{2s-3}}$$ whp.
\end{theorem}

\begin{proof}
The result is immediate when     $p\leq((\ln n)^{\gamma+s-2}/n)^{1/(s-1)}$,   since   in this case  $e(\mathbbmsl{G}(n,p))\ll w_n(\ln n)^{2(\gamma+s-2)}/p^{2s-3}$ whp. If  $p\geq n^{-1/(2s-3)}(\ln n)^\xi$, then   the result follows from Remark \ref{rmrk}. So,  we may assume that
$$\left(\frac{(\ln n)^{\gamma+s-2}}{n}\right)^{\frac{1}{s-1}}<p<n^{-\frac{1}{2s-3}}(\ln n)^\xi.$$
Let $G$ be a graph and let  $X$ be a subset of  $V(G)$ having  two following properties:
\begin{itemize}
\item [(i)] For  any   vertex $u\in V(G)\setminus X$,     $G[N_G(u)\cap X]$  contains the  $(s-2)$-th power of a Hamilton path.
\item [(ii)] For every  two distinct vertices $v,w\in V(G)\setminus X$, $N_G(\{v, w\})\cap X$  contains a clique of size $s-2$.
\end{itemize}
According to  the property (i), for each    $u\in V(G)\setminus X$, we may   fix a     Hamilton  path  $P_u$  of  $G[N_G(u)\cap X]$ such that $P_u^{s-2}\subset G[N_G(u)\cap X]$.
Let	$H$ be a spanning subgraph of $G$ containing all edges with both endpoints in $X$ and also  all $(s-2)|V(G)\setminus X|$   edges in
$$\left\{uv  \, \left| \,
\begin{array}{ll}
\text{$u\in V(G)\setminus X$   and $v$    is one of  the    $s-2$ initial}\\
\text{vertices of $P_u$ starting from  a beginning   vertex}
\end{array}
\right.\right\}.$$
Let us show that $H$ is  a weakly $K_s$-saturated subgraph of  $G$.
Note that all edges with both endpoints in $X$ are initially  activated.
For each $u\in V(G)\setminus X$,
starting from   a beginning   vertex,
denote the vertices of $G[N_G(u)\cap X]$  going in the natural order induced
by $P_u$ by
$v_1, \ldots, v_{x_u}$, where $x_u=|N_G(u)\cap X|$.  The edges  $uv_1,\ldots,uv_{s-2}$ are initially    activated and so
we may active  the edges $uv_{s-1},  \ldots, uv_{x_v}$ one by one,  since   every such edge belongs to a $K_s$  containing the previous $s-2$ vertices of $P_u^{s-2}$.
Hence,  all edges of $G$ going out   $X$   are now activated. Let $v,w\in V(G)\setminus X$ be adjacent in $G$. The property (ii) ensures that the edge  $vw$   belongs to a   $K_s$ whose  other
edges    are already
activated and so the edge $vw$ can be activated as well. This shows that  $H$ is  a weakly $K_s$-saturated subgraph of  $G$ which implies that
$\mathrm{wsat}(G, K_s)\leq e(G[X])+(s-2)|V(G)\setminus X|$.

Let $G\sim\mathbbmsl{G}(n, p)$ and  let  $X=[\![m]\!]$, where $m=\lfloor(\ln n)^{\gamma+s-2}\sqrt{w_n}/p^{s-1}\rfloor$.   Since   $e(G[X])\leq m^2p$ whp, it remains to show that  the properties   (i) and (ii) hold  whp  for $G$ and our choice of $X$.

By Theorem \ref{th:Hamilton_powers} and applying the same coupling technique  as in the proof of Lemma \ref{12},  we   conclude    that,   if    $\ell$ is sufficiently large and     $p\gg\ell^{-1/(s-2)}(\ln n)^{\gamma/(s-2)+1}$, then $\mathbbmsl{G}(n,p)[\ell]$ contains the $(s-2)$-th power of a Hamilton path with probability at least $1-o(\tfrac{1}{n})$. Given a vertex $v\in[\![n]\!]\setminus X$, with probability $o(\tfrac{1}{n})$, $v$ has  at most  $\lfloor(\ln n)^{\gamma+s-2}\sqrt[4]{w_n}/p^{s-2}\rfloor$ neighbors in $X$ using  Theorem \ref{chernoff}. This shows that  the property     (i) holds whp.

Finally, using   Theorem \ref{Janson}, the probability that  there exist distinct  vertices  $v,w\in[\![n]\!]\setminus X$ such that $N_G(\{v, w\})\cap [\![\widetilde m]\!]$ does not   contain  a clique of size $s-2$   is at most
$${n\choose 2}\exp\left(-{\widetilde m\choose s-2}p^{{s\choose 2}-1}+\sum_{t=1}^{s-3}{\widetilde m\choose s-2}{s-2\choose t}{\widetilde m-s+2\choose s-2-t}p^{(s-2)(s+1)-2t-{t\choose 2}}\right),$$ where $\widetilde{m}=\lfloor (\ln n)/p^{(s+1)/2}\rfloor$. The above term is at most
$$\exp\left(2\ln n-\frac{1}{(s-2)!}\widetilde m^{s-2}p^{\frac{(s-2)(s+1)}{2}}\big(1-o(1)\big)\right)$$ which tends to $0$ as $n\to\infty$. Since $\widetilde{m}<m$,  the property  (ii) occurs  whp.
\end{proof}

\section*{Acknowledgments}

Ali Mohammadian  is supported by the    Natural Science Foundation of Anhui Province  with  grant identifier 2008085MA03 and by the National Natural Science Foundation of China with  grant number 12171002.
The works   of Ali Mohammadian  and  Behruz Tayfeh-Rezaie  are  based upon research funded by Iran  National  Science Foundation   under project number  99003814.
Maksim Zhukovskii gratefully acknowledge the financial support from the Ministry of Educational and Science of the Russian Federation in the framework of MegaGrant number  075-15-2019-1926.
The authors wish  to thank the anonymous referees for their valuable comments and  suggestions which helped to improve the presentation of this article.

\appendix\section{Proof of Lemma \ref{edge_in_clique}}\label{appen}

It could be that  Lemma \ref{edge_in_clique}  is a known result in the  literature.  Since we  could not  find any reference, we include its proof here.

\begin{proof}[Proof of Lemma \ref{edge_in_clique}]
By Theorem 2 of  \cite{SP90},  there exists a constant $C>0$ such that, if    $p\geq Cq_s$,   then   every pair of vertices of  $\mathbbmsl{G}(n,p)$     have   whp  $s-2$ common neighbors that induce a clique. Thus, we may assume that $w_n\leq C\ln n$.

Let $N_s$ be the number of edges of  $\mathbbmsl{G}(n,p)$   that do not belong to any $K_s$.
For every two distinct vertices $u,v\in[\![n]\!]$ and any  set   $W\subseteq[\![n]\!]\setminus\{u,v\}$ of size $s-2$, consider the  event $K[W]$ saying that $W$ is a clique in $N_{\mathbbmsl{G}(n,p)}(\{u, v\})$.  Let $\mu(u, v)$ counts the number of sets $W$ as above such that $K[W]$ occurs. We have
\begin{equation}\label{expectation-mu}
\mathbbmsl{E}[N_s]={n\choose 2}p\mathbbmsl{P}\big[\mu(1,2)=0\big].
\end{equation}
If $s=3$, then $\mathbbmsl{P}[\mu(1,2)=0]=(1-p^2)^{n-2}$. For $s\geq 4$, easy  calculations   show that
$$\lambda=\sum_{W\in{[\![n]\!]\setminus\{u, v\}\choose s-2}}\mathbbmsl{P}\big[K[W]\big]={n-2\choose s-2}p^{\frac{(s-2)(s+1)}{2}}$$  and
$$\mathnormal{\Delta}=\mathop{\mathop{\sum_{W_1, W_2\in{[\![n]\!]\setminus\{u,v\}\choose s-2}}}_{W_1\neq W_2}}_{W_1\cap W_2\neq\varnothing}\mathbbmsl{P}\big[K[W_1]\cap K[W_2]\big]$$ is equal to
$$\sum_ {\ell=1}^{s-3}{n-2\choose s-2}{s-2\choose\ell}{n-s\choose s-2-\ell}p^{(s-2)(s-3)-\frac{\ell(\ell-1)}{2}}p^{2(2(s-2)-\ell)}.$$
It follows from  Theorem \ref{FKG} and Theorem \ref{Janson} that
$$\exp\left(-\lambda\left(1+O\left(p^{\frac{(s-2)(s+1)}{2}}\right)\right)\right)\leq\mathbbmsl{P}\big[\mu(1,2)=0\big]\leq
\exp\left(-\lambda+\frac{\mathnormal{\Delta}}{2}\right).$$
According to  the considered expression for  $p(n)$ in  Lemma \ref{edge_in_clique}, we have $\mathnormal{\Delta}=o(\lambda n^{-1/(s+1)}\ln n)$. Hence, since $\lambda=O(\ln n)$,  we deduce  that
\begin{align}
\mathbbmsl{P}\big[\mu(1,2)=0\big]&=\exp\left(-\lambda\left(1+o\left(n^{-\frac{1}{s+1}}\ln n\right)\right)\right)\nonumber\\&=\exp(-\lambda)\exp\left(o\left(n^{-\frac{1}{s+1}}(\ln n)^2\right)\right)\nonumber\\ &=\exp(-\lambda)\left(1+o\left(n^{-\frac{1}{s+1}}(\ln n)^2\right)\right).\label{mu}
\end{align}
If $w_n\to\infty$, then
\begin{align*}
\lambda&=\frac{n^{s-2}}{(s-2)!}\left(2\left(1-\frac{1}{s+1}\right)(s-2)!\right)n^{-(s-2)}\ln n\left(1+\frac{\ln\ln n}{s(s-2)\ln n}+\omega\left(\frac{1}{\ln n}\right)\right)
\\&=2\left(1-\frac{1}{s+1}\right)\ln n+\frac{2\ln\ln n}{(s-2)(s+1)}+\omega(1).
\end{align*}
Therefore, for any $s\geq3$,
\begin{align*}
\mathbbmsl{E}[N_s]&\leq{n\choose 2}p\exp\left(-\lambda+\frac{\mathnormal{\Delta}}{2}\right)\\&\leq\exp\big(2\ln n+\ln p-\lambda+O(1)\big)\\
&=\exp\left(2\ln n-\frac{2\ln n}{s+1}+\frac{2\ln\ln n}{(s-2)(s+1)}-2\left(1-\frac{1}{s+1}\right)\ln n-\frac{2\ln\ln n}{(s-2)(s+1)}-\omega(1)\right),
\end{align*}
which goes to $0$ as $n\to\infty$ and thus  Part (i)  of Lemma \ref{edge_in_clique} follows by Corollary \ref{First-Moment}.

In order  to  prove  Part (ii)  of Lemma \ref{edge_in_clique}, note first that,  if $w_n\to-\infty$, then
$$\mathbbmsl{E}[N_s]\geq{n\choose 2}p\exp\left(-\lambda\left(1+O\left(n^{-(s-2)}\ln n\right)\right)\right)\geq\exp\big(2\ln n+\ln p-\lambda+o(1)\big)\to\infty$$ as $n\to\infty$, meaning that  $\mathbbmsl{E}[N_s]=o(\mathbbmsl{E}[N_s]^2)$.
Below, we   apply   Theorem \ref{Janson} to   estimate $\mathbbmsl{P}[\mu(1,2)=\mu(3,4)=0]$ and $\mathbbmsl{P}[\mu(1,2)=\mu(2,3)=0]$.
For each set $W\subseteq[\![n]\!]\setminus[\![4]\!]$ of size $s-2$, consider the  event $K'[W]$ saying that   $W$ is a clique in either $N_{\mathbbmsl{G}(n,p)}(\{1, 2\})$ or   $N_{\mathbbmsl{G}(n,p)}(\{3, 4\})$.   Let $\mu'$ counts the number of sets $W$ as above such that $K'[W]$ happens. We have
\begin{equation}
\mathbbmsl{P}\big[\mu(1,2)=\mu(3,4)=0\big]\leq\mathbbmsl{P}[\mu'=0].\label{mu-one}
\end{equation}
If $s=3$, then $\mathbbmsl{P}[\mu'=0]=(1-p^2)^{2(n-4)}$. For $s\geq 4$, easy  computations  show that
$$\lambda'=\sum_{W'\in{[\![n]\!]\setminus[\![4]\!]\choose s-2}}\mathbbmsl{P}\big[K'[W]\big]={n-4\choose s-2}p^{\frac{(s-2)(s+1)}{2}}\left(2-p^{2(s-2)}\right)$$   and
$$\mathnormal{\Delta}'=\mathop{\mathop{\sum_{W_1,W_2\in{[\![n]\!]\setminus[\![4]\!]\choose s-2}}}_{W_1\neq W_2}}_{W_1\cap W_2\neq\varnothing}\mathbbmsl{P}\big[K'[W_1]\cap K'[W_2]\big]$$ is equal to $$\sum_{\ell=1}^{s-3}{n-4\choose s-2}{s-2\choose\ell}{n-s-2\choose s-2-\ell}
p^{(s-2)(s-3)-\frac{(\ell-1)\ell}{2}}2p^{2\ell}\left(2p^{2(s-2)-2\ell}\right)^2\big(1+O(p)\big).$$
We have $\lambda'=2\lambda(1+O(\tfrac{1}{n}))$ and $\mathnormal{\Delta}'=8\mathnormal{\Delta}(1+O(q_s))$.
So, Theorem \ref{Janson} implies that
\begin{align}
\mathbbmsl{P}[\mu'=0]&\leq\exp\left(-\lambda'+\frac{\mathnormal{\Delta}'}{2}\right)\nonumber\\&=\exp\left(-2\lambda\left(1+O\left(\frac{1}{n}\right)\right)+4\mathnormal{\Delta}\big(1+O(q_s)\big)\right)\nonumber\\&=
\exp\left(-2\lambda\left(1+o\left(n^{-\frac{1}{s+1}}\ln n\right)\right)\right)\nonumber\\&=\exp(-2\lambda)\left(1+o\left(n^{-\frac{1}{s+1}}(\ln  n)^2\right)\right).\label{mu-one-from-above}
\end{align}
In the same way, let us estimate $\mathbbmsl{P}[\mu(1,2]=\mu(2,3)=0]$.
For each  set $W\subseteq[\![n]\!]\setminus[\![3]\!]$ of size $s-2$, consider the  event $K''[W]$ saying that   $W$ is a clique in either $N_{\mathbbmsl{G}(n,p)}(\{1, 2\})$ or   $N_{\mathbbmsl{G}(n,p)}(\{2, 3\})$. Let $\mu''$ counts the number of sets $W$ as above such that $K''[W]$ happens. We have
\begin{equation}
\mathbbmsl{P}\big[\mu(1,2)=\mu(2,3)=0\big]\leq\mathbbmsl{P}[\mu''=0].\label{mu-two}
\end{equation}
If $s=3$, then $\mathbbmsl{P}[\mu''=0]=(1-p+p(1-p)^2)^{n-3}$. For $s\geq 4$, easy  calculations   show that
$$\lambda''=\sum_{W'\in{[\![n]\!]\setminus[\![3]\!]\choose s-2}}\mathbbmsl{P}\big[K''[W]\big]={n-3\choose s-2}p^{\frac{(s-2)(s+1)}{2}}\left(2-p^{s-2}\right)$$ and
$$\mathnormal{\Delta}''=\mathop{\mathop{\sum_{W_1,W_2\in{[\![n]\!]\setminus[\![3]\!]\choose s-2}}}_{W_1\neq W_2}}_{W_1\cap W_2\neq\varnothing}\mathbbmsl{P}\big[K''[W_1]\cap K''[W_2]\big]$$ is equal to
$$\sum_{\ell=1}^{s-3}{n-3\choose s-2}{s-2\choose\ell}{n-s-1\choose s-2-\ell}p^{(s-2)(s-3)-\frac{(\ell-1)\ell}{2}}p^{2(s-2)-\ell}8p^{\ell}p^{2(s-2)-2\ell}\big(1+O(p)\big).$$
We have $\lambda''=2\lambda(1+o(\tfrac{1}{\sqrt{n}}))$ and $\mathnormal{\Delta}''=8\mathnormal{\Delta}(1+O(q_s))$.
Then, Theorem \ref{Janson} implies that
\begin{align}
\mathbbmsl{P}[\mu''=0]&\leq\exp\left(-\lambda''+\frac{\mathnormal{\Delta}''}{2}\right)\nonumber\\&=\exp\left(-2\lambda\left(1+O\left(\frac{1}{\sqrt{n}}\right)\right)+4\mathnormal{\Delta}\big(1+O(q_s)\big)\right)
\nonumber\\&=\exp\left(-2\lambda\left(1+o\left(n^{-\frac{1}{s+1}}\ln n\right)\right)\right).\label{mu-two-from-above}
\end{align}
By combining the relations \eqref{expectation-mu}--\eqref{mu-two-from-above}, we derive  for each  $s\geq3$  that
\begin{align*}
\mathbbmsl{Var}[N_s]&=\mathbbmsl{E}[N_s^2]-\mathbbmsl{E}[N_s]^2\\&=\mathbbmsl{E}[N_s]+{n\choose 2}{n-2\choose 2}p^2\mathbbmsl{P}\big[\mu(1,2)=\mu(3,4)=0\big]\\&+
n(n-1)(n-2)p^2\mathbbmsl{P}\big[\mu(1,2)=\mu(2,3)=0\big]-\left({n\choose 2}p\mathbbmsl{P}\big[\mu(1,2)=0\big]\right)^2\\&\leq
\mathbbmsl{E}[N_s]+{n\choose 2}^2p^2\exp(-2\lambda)\left(1+o\left(n^{-\frac{1}{s+1}}(\ln n)^2\right)\right)\\&+
n^3p^2\exp(-2\lambda)\left(1+o\left(n^{-\frac{1}{s+1}}(\ln n)^2\right)\right) - {n\choose 2}^2p^2\exp(-2\lambda)\left(1+o\left(n^{-\frac{1}{s+1}}(\ln n)^2\right)\right)\\&=o\big(\mathbbmsl{E}[N_s]^2\big).
\end{align*}
Now, Part (ii)  of Lemma \ref{edge_in_clique}  follows from Corollary   \ref{Second-Moment}.
\end{proof}

\end{document}